\documentclass[11pt]{amsart}
\usepackage{amsmath,amsfonts}
\usepackage{amssymb}
\usepackage{amscd}
\usepackage{amsthm}
\usepackage{yhmath}
\usepackage{subfigure}

\usepackage[all]{xy}
\usepackage{color}

\numberwithin{equation}{section}

\newtheorem{theorem}[equation]{Theorem}

\theoremstyle{remark}

\theoremstyle{definition}
\newtheorem{definition}[equation]{Definition}

\def\XXint#1#2#3{{\setbox0=\hbox{$#1{#2#3}{\int}$}
	\vcenter{\hbox{$#2#3$}}\kern-.5\wd0}}

\newcommand{\R}{\mathbb R}

\newcommand{\G}{{\mathbb G}}

\newcommand{\norm}[1]{\left\Vert#1\right\Vert}

\begin{document}

\title{A metric characterization of Carnot groups}

\author{Enrico Le Donne}
\address{Department of Mathematics and Statistics \\
         P.O. Box 35 (MaD) \\
         FI-40014 
         University of Jyv\"askyl\"a \\
         Finland}
\email{enrico.e.ledonne@jyu.fi}

\thanks{The author thanks IPAM and all of the people involved in the program `Interactions Between Analysis and Geometry', 
during which  there was the 
 opportunity 
 of discussing these results.
 }

\renewcommand{\subjclassname}{%
 \textup{2010} Mathematics Subject Classification}
\date{August 25, 2014}

\begin{abstract} 
We give a short axiomatic  introduction to Carnot groups and their subRiemannian and
 subFinsler geometry.
We  explain how such spaces can be metrically described as exactly
 those proper geodesic spaces that admit dilations and are
 isometrically homogeneous.
 \end{abstract}

\maketitle

\section{Introduction}
 
 Carnot groups are special  objects that play important roles  in several mathematical fields.
 For example, they appear in Algebra and Geometric Group Theory; in particular, in the  theory of nilpotent groups, and 
 in  Geometric Analysis and PDE,  as models for  
subelliptic operators.

Briefly speaking, 
 a Carnot group, or better a subFinsler Carnot group, is a nilpotent stratified Lie group equipped with a left-invariant subFinsler distance with the first stratum as horizontal distribution.
 The purpose of this short paper is to give a more axiomatic presentation of Carnot groups from the view point of Metric Geometry. In fact, we shall see that Carnot groups are the only locally compact and geodesic metric spaces that are isometrically homogeneous and self-similar.
 Such a result follows the spirit of Gromov's approach of `seeing {C}arnot-{C}arath\'eodory spaces  from within', \cite{Gromov1}.

Let us recall and make explicit the above definitions.
A topological space $X$ is called {\em locally compact} if every point of the space has a compact neighborhood. 
A metric space is {\em geodesic} if, for all $p,q\in X$, there exists an isometric embedding $\iota:[0,T]\to X$ with $T\geq 0 $ such that 
$\iota(0)=p$ and $\iota(T)=q$.
We say that a metric space $X$ is {\em isometrically homogeneous} if its group of isometries acts on the space transitively. Explicitly, this means that, for all $p,q\in X$, there exists a distance-preserving homeomorphism $f:X\to X$ such that $f(p)=q$.
In this paper, we  say that  a metric space $X$ is {\em  self-similar} if it admits a {\em dilation}, i.e., there exists $\lambda>1$ and a homeomorphism  $f:X\to X$ such that $d(f(p), f(q) )=\lambda d(p,q) ,\text{  for all }p,q\in X$.
 \begin{samepage}
 \begin{theorem}\label{characterization_Carnot}
The subFinsler Carnot groups are the only metric spaces that are  
\begin{enumerate}
\item  locally compact, 
\item  geodesic, 
  \item isometrically homogeneous, and
   \item  self-similar (i.e.,  admitting a dilation).
\end{enumerate}
 \end{theorem}
 \end{samepage}
 
 
 Theorem \ref{characterization_Carnot} provides a new equivalent definition of Carnot groups. 
 Obviously, (1) can be slightly strengthened assuming that the space is {\em boundedly compact} (the term {\em proper} is also used), i.e., closed balls are compact. 
 Let us now recall  what is the 
 traditional definition of a Carnot group.
  \begin{definition}[SubFinsler Carnot groups]
 Let $\G$ be a connected and simply-connected Lie group.
 Denote by ${\rm Lie}(\G)$ the Lie algebra of $\G$, seen as the set of   tangent vectors at the identity element with the bracket operation induced by the Lie bracket of left-invariant vector fields.
The group $\G$   is said to be a   {\em stratified group} 
if 
there exist subspaces $V_j\subseteq {\rm Lie}(\G)$ such that
\begin{equation}
\label{stratification}
{\rm Lie}(\G)= V_1\oplus \cdots\oplus V_s, \quad \text{ with }\, [V_j, V_1] = V_{j+1}, \quad\text{ for } 1\leq j\leq s,
\end{equation} where   $V_{s+1}= \{0\}$.
If $V_s\neq \{0\}$, we say that \eqref{stratification} is a {\em stratification of step} $s$.
Note that such a group $\G$ is nilpotent and that the {\em first stratum} $V_1$ generates the whole Lie algebra.

The space $V_1$ defines a left-invariant subbundle $\Delta$ of the tangent bundle of $\G$, which is named {\em horizontal distribution}. Namely,
$$\Delta_p := (L_p)_*V_1,$$
where $L_p$ is the left-multiplication by $p \in \G$ and $F_*$ denotes the differential of a diffeomorphism $F$.
Given a norm  $ \norm{\cdot} $ on $V_1$, we can extend it left-invariantly on $\Delta $ as
$$ \norm{v} := \norm{(L_p^{-1})_*  v} , \qquad \forall v\in \Delta_p  . $$ 
For an absolutely continuous  curve $\gamma:[0,1]\to \G$ that has the property that almost everywhere  $\dot\gamma\in \Delta$, which is called a {\em horizontal curve}, we set
$$
    {\rm Length}_ {\norm{\cdot}}(\gamma) := 
    \int_0^1 \norm{\dot\gamma(t)}\,\mathrm d t
    .
$$
  The associated \textit{subFinsler}  (or 
 \textit{Finsler-Carnot-Carath\'eodory}) distance 
  between two points $p,q\in \G$  is 
  defined\footnote{As shown in \cite[Theorem 1.2]{LeDonne1}, in
\eqref{dist_CC} one can replace the smoothness condition on $\gamma$ with just  horizontality. 
} as
\begin{equation}\label{dist_CC}
d(p,q):=\inf\{{\rm Length}_ {\norm{\cdot}}(\gamma)\;|\;\gamma\in C^\infty([0,1];\G), \gamma(0)= p , \gamma(1)= q, \dot \gamma\in \Delta\}.
\end{equation}
We call the metric space
$(\G,d)$ a
{\em subFinsler Carnot group}.

If the chosen norm   comes from a scalar product, then the associated distance is called {\it Carnot-Carath\'eodory} or {\it subRiemannian}. 
If this is the case,  we call $(G,d)$ a {\em subRiemannian Carnot group}.  
The term Carnot group is often used in the literature. This ambiguous terminology, which will also be used here,
  might denote either a subRiemannian or, more generally, a subFinsler Carnot group. Some other authors  use the term Carnot group even as a synonym for stratified group.
\end{definition}

The topology of the metric space
$(\G,d)$ is the same of the one of $\G$ as manifold. This last fact is shown using the fact that $V_1$ generates the Lie algebra, and it holds more generally; see Chow-Rashevsky
Theorem
 \cite{Montgomery}. Hence, since a Carnot group is locally compact,  complete, and the distance is defined as a length distance,  it is a geodesic space.
  Obviously, a Carnot group is isometrically homogeneous: since the whole construction of the distance is left-invariant,  left translations are isometries (and act transitively).
   Regarding  self-similarity, Carnot groups  admit  dilations for any factor $\lambda$, as we now recall. This property will be a consequence of the fact that the group is stratified (and simply connected).
  At the Lie algebra level, for each $\lambda \in \R$, the dilation $\tilde\delta_\lambda$  is defined linearly by setting $\tilde\delta_\lambda(X):=\lambda^jX$, for every $X\in V_{j}$  and every $j=1,\dots,s$.  For all $\lambda\neq0$,     
the map $\tilde\delta_\lambda$ is a Lie algebra isomorphism.
From the theory of    simply connected   Lie groups, we have that the dilation induces a unique isomorphism of the group, which we denote by $\delta_\lambda$. Namely, there exists a map $\delta_\lambda :\G\to \G$ such that $(\delta_\lambda)_*=\tilde\delta_\lambda$. 
Hence, the map $\delta_\lambda$ stretches the length of horizontal curves exactly by $\lambda$. Therefore, it dilates the distance by $\lambda$.
Observe that, since  $\G$ is nilpotent and simply connected, the exponential map is a diffeomorphism; thus the dilations   $\delta_\lambda$ can be equivalently defined as
 $\delta_\lambda(p)=\exp \circ \, \tilde\delta_\lambda \circ \exp^{-1}(p)$, for all $\lambda\in \R$ and  $p \in \G$. 
 
 We would like to remark a few things. Note that a Carnot group admits dilations for {\em all} factors $\lambda$.
Nonetheless, we 
claim that having the existence of {\em  one} dilation is enough, together with the other assumptions, to conclude the existence of  dilations of   any factor.

We point out that each of the four conditions in Theorem \ref{characterization_Carnot} is necessary for the validity of the result. Indeed, 
let us mention examples of spaces that satisfy three out of the four conditions but are not Carnot groups:
any infinite dimensional Banach space;
any snowflake of a Carnot group, e.g., $(\R, \sqrt{\norm{\cdot}})$;
many cones such as the usual Euclidean cone of cone angle in $(0,2\pi)$ or the union of two spaces such as
$\{(x,y)\in \R^2: xy\geq0\}$;
any compact homogeneous space such as $\mathbb S^1$.

This paper is not the first one focusing on metric characterizations of Carnot groups. Other papers in this context are
\cite{LeDonne6}, \cite{Buliga11}, \cite{Freeman12} (which is based on \cite{LeDonne1}), and \cite{Buyalo-Schroeder_characeterization}.
 
 \subsection*{Examples of Carnot groups}
 Here are very basic samples of Carnot groups. 

1. The first examples are the Euclidean spaces and the finite dimensional Banach spaces. Here the group structure is Abelian, hence the stratification has step one, i.e., there is only one stratum, which coincides with the whole Lie algebra.

2. The next very popular Carnot group is the Heisenberg group, which we denote by $\mathbb H$.
The space $\mathbb H$ is $\R^3$ equipped with a distance as follows.
Fix a norm  $ \norm{\cdot} $ on $\R^2$.
Then define a distance between two points $p$ and $q$ of $\R^3$ as
\begin{equation*}\label{dist_CC-Heis}
d(p,q):=\inf\left\{   \left.
    \int_0^1 \norm{(\dot\gamma_1(t),\dot\gamma_2(t))  }\mathrm d t
    \;\right|\;\gamma\in C^\infty([0,1];\R^3), 
    \begin{array}{ccc}
	 \gamma(0)= p, \\
	\gamma(1)= q, 
	\end{array}
   \dot \gamma_3 =\frac{1}{2}( \gamma_1 \dot \gamma_2-\gamma_2 \dot \gamma_1)
   \right \}.
\end{equation*}
The group structure, for which this distance is left-invariant, is given by 
$$(\bar x,\bar y,\bar z) \cdot (x,y,z):= \left(\bar x+x,\bar y+y,\bar z+z+\frac{1}{2} (\bar x y-\bar y x)\right).          $$
A dilation by factor $\lambda$ is given by the map
$(x,y,z) \longmapsto (\lambda x,\lambda y,\lambda^2 z) $.
For this group we have the step-two stratification 
$$V_1:=\text{span}\left\{\partial_1-\dfrac{y}{2}\partial_3,\, \partial_2+\dfrac{x}{2}\partial_3 \right\} \qquad \text{ and } \quad\qquad V_2:=\text{span}\{\partial_3\}.$$

 \section{Proof of the characterization}
 The proof of Theorem \ref{characterization_Carnot} is an easy consequence of three hard theorems. 
 We present now these theorems, before giving the proof.
	
The first theorem is well-known in the theory of locally compact groups. It is a consequence of a deep result of Dean Montgomery and Leo Zippin,
\cite[Corollary on page 243, Section 6.3]{mz}, together with the work \cite{Gleason} of Andrew Gleason.
An explicit proof can be found in Cornelia Drutu and Michael Kapovich's lecture notes, \cite[Chapter 14]{Drutu-Kapovich}.

\begin{theorem}[Gleason-Montgomery-Zippin]\label{Montgomery-Zippin}
Let $X$  be a metric space that is connected,  locally connected,  locally compact and has finite topological dimension.
Assume that the isometry group Isom$(X)$ of $X$ acts transitively on $X$.
Then Isom$(X)$  has the  structure of a Lie group with finitely many
connected components,  and $X$ has the structure of an analytic manifold.
\end{theorem}
Notice that  
 an isometrically homogeneous space that is locally compact is complete.

Successively, Berestovskii's work  \cite[Theorem 2]{b} clarified what  the possible isometrically homogeneous distances on manifolds that are also geodesic are. They are subFinsler metrics.
\begin{theorem}[Berestovskii]\label{Berestovskii}
Under the same assumptions of  Theorem \ref{Montgomery-Zippin}, if in addition the distance is geodesic, then the distance is a subFinsler metric, i.e.,  the metric space $X$ is a homogeneous Lie space $G/H$ and   there is a $G$-invariant subbundle $\Delta$ on the manifold $G/H$ and a $G$-invariant norm on $\Delta$ such that the distance is given by the same formula \eqref{dist_CC}. 
\end{theorem}
 
Tangents, in the Gromov-Hausdorff sense, of subFinsler manifolds have been studied. In particular, recall that a bundle is {\em equiregular} if  the spaces of iterated brackets form subbundles as well. For  more detailed terminology and for the proof of the next result we refer to  \cite{bellaiche, 
Margulis-Mostow, Margulis-Mostow2}.
\begin{theorem}[Mitchell]\label{Mitchell}
The metric tangents of an equiregular subFinsler manifold are subFinsler Carnot groups.
\end{theorem}

\subsection*{Proof of Theorem \ref{characterization_Carnot}}
Let us verify that we can use Theorem \ref{Montgomery-Zippin}. A geodesic metric space is obviously   connected and  locally connected.
Regarding finite dimensionality, we claim that a locally compact, self-similar, isometrically homogeneous space $X$ is doubling. Namely, there exists a constant $C>0$ such that any ball of radius $r>0$ in $X$ can be covered with less than $C$ balls of radius $r/2$. Since $X$ is locally compact, there exists a ball $B(x_0,r_0)$ that is compact.
Let $\lambda>1$ be the factor of the dilation. Hence,  the balls $B(x_0, s r_0)$ with $s\in [1,\lambda]$ form a compact family of compact balls. 
 Hence, there exists a constant $C>1$ such that
 each  ball $B(x_0, s r_0)$ can be covered   with less than $C$ balls of radius $s r_0/2$.  By self-similarity and homogeneity, any other ball can be  covered with less than $C$ balls of half radius. Doubling metric spaces  have finite Hausdorff dimension and hence finite  topological dimension. 
Therefore, by Theorem \ref{Montgomery-Zippin} the isometry group $G$ is a Lie group.

Since the distance is  geodesic, Theorem \ref{Berestovskii} implies that our metric space is a subFinsler homogeneous manifold  $G/H$. Since the subFinsler structure is $G$ invariant, in particular it is equiregular. Hence, on the one hand, because of 
Theorem \ref{Mitchell}
the tangents of our metric space are subFinsler Carnot groups.
On the other hand,  the space admits a dilation; hence, iterating the dilation, we have that there exists a metric tangent of the metric space that is isometric to our original space.
Then the space is a subFinsler Carnot group.
\qed

After completion of this manuscript, Valerii Berestovskii has informed the author that a statement similar to Theorem \ref{characterization_Carnot} can be found in his work \cite{Berestovskii_2004}.

   \bibliography{general_bibliography-copy}
\bibliographystyle{amsalpha}

\end{document}